# Symmetry operators and generation of symmetry transformations of partial differential equations

**Costas J. Papachristou**

Department of Physical Sciences, Hellenic Naval Academy, Piraeus, Greece
E-mail: cjpapachristou@gmail.com

**Abstract.** The study of symmetries of partial differential equations (PDEs) has been traditionally treated as a geometrical problem. Although geometrical methods have been proven effective with regard to finding infinitesimal symmetry transformations, they present certain conceptual difficulties in the case of matrix-valued PDEs; for example, the usual differential-operator representation of the symmetry-generating vector fields is not possible in this case. In this article an algebraic approach to the symmetry problem of PDEs – both scalar and matrix-valued – is described, based on abstract operators (characteristic derivatives) that admit a standard differential-operator representation in the case of scalar-valued PDEs. A number of examples are given.

## 1. Introduction

The symmetry properties of partial differential equations (PDEs) have been traditionally treated as a geometrical problem in the jet space of the independent and the dependent variables (including a sufficient number of partial derivatives of the latter variables with respect to the former ones). Two more or less equivalent approaches have been adopted: (*a*) invariance of a system of PDEs under infinitesimal transformations generated by corresponding vector fields in the jet space [1]; (*b*) invariance of a differential ideal of differential forms representing the system of PDEs, under the Lie derivative with respect to the vector fields representing the symmetry transformations [2-6].

Although effective with regard to calculating symmetries, these geometrical approaches suffer from a certain drawback of conceptual nature when it comes to *matrix-valued* PDEs. The problem is related to the inevitably mixed nature of the variables in the jet space (scalar independent variables versus matrix-valued dependent ones) and the need for a differential-operator representation of the symmetry vector fields. How does one define differentiation with respect to matrix-valued variables? Moreover, how does one calculate the Lie bracket of two differential operators in which some (or all) of the variables, as well as the coefficients of partial derivatives with respect to these variables, are matrices?

Although these difficulties were handled geometrically in a somewhat *ad hoc* way in [4-6], it was eventually realized that an alternative, *purely algebraic* approach to the symmetry problem would be more appropriate in the case of matrix-valued PDEs. Such an approach was presented in [7], although in practice it had already been applied in particular problems [8-10] before a full-blown formal theoretical framework was developed.

This algebraic framework for symmetries of PDEs is described in this article. In Sections 2 and 3 we introduce the concept of *characteristic derivatives* – an abstract generalization of the geometrical notion of vector fields in differential-operator form – and we demonstrate the Lie-algebraic character of the set of these derivatives.

The general symmetry problem for both scalar and matrix-valued PDEs is presented in Sec. 4, and the Lie-algebraic property of symmetries of a PDE is proven in Sec. 5. In Sec. 6 we discuss the concept of a *recursion operator* [1,8-14] by which an infinite set of symmetries may in principle be produced from any known symmetry. An application of these ideas is made in Sec. 7 by using the chiral field equation as an example.

A symmetry of a PDE amounts to the invariance of this equation under the action of a corresponding characteristic derivative. Given the latter operator, an *infinitesimal* symmetry of the PDE may be defined. Section 8 discusses the use of symmetry operators to construct *finite* one-parameter symmetry transformations of PDEs. As a pedagogical example, a number of point-symmetry transformations for the two-dimensional Laplace equation are derived in Sec. 9. Finally, in Sec. 10 the case of Bäcklund-transformation-related recursion operators [9] is discussed.

To simplify our formalism we will restrict our analysis to the case of a single PDE in one dependent variable. For systems of scalar-valued PDEs in several dependent variables see, e.g., [1].

The material presented in this article is largely based on the original article [7], with additional results contained in Sec. 8-10.

## 2. The fundamental operators

A PDE for the unknown function $u=u(x^1, x^2, \ldots) \equiv u(x^k)$ [where by $(x^k)$ we collectively denote the independent variables $x^1, x^2, \ldots$] is an expression of the form $F[u]=0$, where $F[u] \equiv F(x^k, u, u_k, u_{kl}, \ldots)$ is a function in the *jet space* [1] of the independent variables $(x^k)$, the dependent variable $u$, and the partial derivatives of various orders of $u$ with respect to the $x^k$, which derivatives will be denoted by using subscripts: $u_k$, $u_{kl}$, $u_{klm}$, etc. A *solution* of the PDE is any function $u=\varphi(x^k)$ for which the relation $F[u]=0$ is satisfied.

The dependent variable $u$, as well as all functions $F[u]$ in the jet space, will generally be assumed to be square-matrix-valued of fixed (but otherwise unspecified) matrix dimensions. In particular, we require that, in its most general form, a function $F[u]$ in the jet space is expressible as a finite or an infinite sum of products of alternating $x$-dependent and $u$-dependent terms, of the form

$$F[u] = \sum a(x^k)\, \Pi[u]\, b(x^k)\, \Pi'[u]\, c(x^k) \cdots \qquad (1)$$

where the $a(x^k)$, $b(x^k)$, $c(x^k)$, etc., are (generally) matrix-valued and where the matrices $\Pi[u]$, $\Pi'[u]$, etc., are products of variables $u$, $u_k$, $u_{kl}$, etc., of the "fiber" space (or, more generally, products of powers of these variables). The set of all functions (1) is thus a (generally) non-commutative algebra.

If $u$ is a scalar quantity, a total derivative operator can be defined in the usual differential-operator form

$$D_i = \frac{\partial}{\partial x^i} + u_i \frac{\partial}{\partial u} + u_{ij} \frac{\partial}{\partial u_j} + u_{ijk} \frac{\partial}{\partial u_{jk}} + \cdots \tag{2}$$

where the summation convention over repeated up-and-down indices (such as $j$ and $k$ in this equation) has been adopted and will be used throughout. If, however, $u$ is matrix-valued, the above expression is obviously not valid. A generalization of the definition of the total derivative is thus necessary for matrix-valued PDEs.

*Definition 1.* The *total derivative operator* with respect to the variable $x^i$ is a linear operator $D_i$ acting on functions $F[u]$ of the form (1) in the jet space and having the following properties:

1. On functions $f(x^k)$ in the base space, $D_i f(x^k) = \partial f / \partial x^i \equiv \partial_i f(x^k)$.

2. For $F[u]=u$, $u_i$, $u_{ij}$, etc., we have: $D_i u = u_i$, $D_i u_j = D_j u_i = u_{ij} = u_{ji}$, etc.

3. The operator $D_i$ is a *derivation* on the algebra of all matrix-valued functions of the form (1) in the jet space; i.e., the *Leibniz rule* is satisfied:

$$D_i \left( F[u] G[u] \right) = \left( D_i F[u] \right) G[u] + F[u] D_i G[u] \tag{3}$$

Higher-order total derivatives $D_{ij} = D_i D_j$ may similarly be defined but they no longer possess the derivation property. Given that $\partial_i \partial_j = \partial_j \partial_i$ and that $u_{ij} = u_{ji}$, it follows that $D_i D_j = D_j D_i \Leftrightarrow D_{ij} = D_{ji}$; that is, total derivatives commute. We write: $[D_i, D_j]=0$, where, in general, $[A, B] \equiv AB-BA$ will denote the *commutator* of two operators or two matrices, as the case may be.

If $u^{-1}$ is the inverse of $u$, such that $uu^{-1} = u^{-1} u = 1$, then we can define

$$D_i \left( u^{-1} \right) \equiv - u^{-1} \left( D_i u \right) u^{-1} \tag{4}$$

Moreover, for any functions $A[u]$ and $B[u]$ in the jet space it can be shown that

$$D_i [A, B] = [D_i A, B] + [A, D_i B] \tag{5}$$

As an example, let $(x^1, x^2) \equiv (x, t)$ and let $F[u]=xtu_x^2$, where $u$ is matrix-valued. Writing $F[u]=xtu_x u_x$, we have: $D_t F[u]=xu_x^2 + xt\,(u_{xt} u_x + u_x u_{xt})$.

Let now $Q[u] \equiv Q(x^k, u, u_k, u_{kl}, \ldots)$ be a function in the jet space. We will call this a *characteristic function* (or simply a *characteristic*) of a certain derivative, defined as follows:

*Definition 2.* The *characteristic derivative* with respect to $Q[u]$ is a linear operator $\Delta_Q$ acting on functions $F[u]$ in the jet space and having the following properties:

1. On functions $f(x^k)$ in the base space,

$$\Delta_Q f(x^k) = 0 \tag{6}$$

(that is, $\Delta_Q$ acts only in the fiber space).

2. For $F[u]=u$,

$$\Delta_Q u = Q[u] \tag{7}$$

3. $\Delta_Q$ commutes with total derivatives:

$$\Delta_Q D_i = D_i \Delta_Q \quad \Leftrightarrow \quad [\Delta_Q , D_i ] = 0 \quad (\text{all } i) \tag{8}$$

4. The operator $\Delta_Q$ is a *derivation* on the algebra of all matrix-valued functions of the form (1) in the jet space (the Leibniz rule is satisfied):

$$\Delta_Q \left( F[u]G[u] \right) = \left( \Delta_Q F[u] \right) G[u] + F[u]\, \Delta_Q G[u] \tag{9}$$

*Corollary:* By (7) and (8) we have:

$$\Delta_Q u_i = \Delta_Q D_i u = D_i Q[u] \tag{10}$$

We note that the operator $\Delta_Q$ is a well-defined quantity, in the sense that the action of $\Delta_Q$ on $u$ uniquely determines the action of $\Delta_Q$ on any function $F[u]$ of the form (1) in the jet space. Moreover, since, by assumption, $u$ and $Q[u]$ are matrices of equal dimensions, it follows from (7) that $\Delta_Q$ preserves the matrix character of $u$, as well as of any function $F[u]$ on which this operator acts.

We also remark that we have imposed conditions (6) and (8) having a certain property of symmetries of PDEs in mind; namely, *every* symmetry of a PDE can be represented by a transformation of the dependent variable $u$ alone, i.e., can be expressed as a transformation in the fiber space (see [1], Chap. 5).

The following formulas, analogous to (4) and (5), may be written:

$$\Delta_Q \left( u^{-1} \right) \equiv - u^{-1} \left( \Delta_Q u \right) u^{-1} \tag{11}$$

$$\Delta_Q [A,B] = [\Delta_Q A, B] + [A, \Delta_Q B] \tag{12}$$

As an example, let $(x^1, x^2) \equiv (x, t)$ and let $F[u]=a(x,t)u^2b(x,t)+[u_x , u_t]$ , where $a$, $b$ and $u$ are matrices of equal dimensions. Writing $u^2=uu$ and using properties (7), (9), (10) and (12), we find: $\Delta_Q F[u]=a(x,t)(Qu+uQ)b(x,t)+[D_x Q, u_t]+[u_x , D_t Q]$.

In the case where $u$ is scalar-valued (thus so is $Q[u]$) the characteristic derivative $\Delta_Q$ admits a differential-operator representation of the form

$$\Delta_Q = Q[u] \frac{\partial}{\partial u} + \left( D_i Q[u] \right) \frac{\partial}{\partial u_i} + \left( D_i D_j Q[u] \right) \frac{\partial}{\partial u_{ij}} + \cdots \tag{13}$$

See [1], Chap. 5, for an analytic proof of property (8) in this case.

## 3. The Lie algebra of characteristic derivatives

The characteristic derivatives $\Delta_Q$ acting on functions $F[u]$ of the form (1) in the jet space constitute a *Lie algebra of derivations* on the algebra of the $F[u]$. The proof of this statement is contained in the following three Propositions.

*Proposition 1.* Let $\Delta_Q$ be a characteristic derivative with respect to the characteristic $Q[u]$; i.e., $\Delta_Q u=Q[u]$ [cf. Eq. (7)]. Let $\lambda$ be a constant (real or complex). We define the operator $\lambda\Delta_Q$ by the relation

$$(\lambda\Delta_Q)F[u] \equiv \lambda\,(\Delta_Q F[u])\,.$$

Then, $\lambda\Delta_Q$ is a characteristic derivative with characteristic $\lambda Q[u]$. That is,

$$\lambda\,\Delta_Q = \Delta_{\lambda Q} \tag{14}$$

*Proof.* (*a*) The operator $\lambda\Delta_Q$ is linear, since so is $\Delta_Q$ .

(*b*) For $F[u]=u$, $(\lambda\Delta_Q)u= \lambda(\Delta_Q u)=\lambda Q[u]$ .

(*c*) $\lambda\Delta_Q$ commutes with total derivatives $D_i$ , since so does $\Delta_Q$ .

(*d*) Given that $\Delta_Q$ satisfies the Leibniz rule (9), it is easily shown that so does $\lambda\Delta_Q$ .

*Comment:* Condition (*c*) would not be satisfied if we allowed $\lambda$ to be a function of the $x^k$ instead of being a constant, since $\lambda(x^k)$ generally does not commute with the $D_i$. Therefore, relation (14) is not valid for a non-constant $\lambda$.

*Proposition 2.* Let $\Delta_1$ and $\Delta_2$ be characteristic derivatives with respect to the characteristics $Q_1[u]$ and $Q_2[u]$, respectively; i.e., $\Delta_1 u=Q_1[u]$, $\Delta_2 u=Q_2[u]$. We define the operator $\Delta_1+\Delta_2$ by

$$(\Delta_1+\Delta_2)\,F[u] \equiv \Delta_1\,F[u] + \Delta_2\,F[u]\,.$$

Then, $\Delta_1+\Delta_2$ is a characteristic derivative with characteristic $Q_1[u]+Q_2[u]$. That is,

$$\Delta_1 + \Delta_2 = \Delta_Q \quad \text{with} \quad Q[u] = Q_1[u] + Q_2[u] \tag{15}$$

*Proof.* (*a*) The operator $\Delta_1+\Delta_2$ is linear, as a sum of linear operators.

(*b*) For $F[u]=u$, $(\Delta_1+\Delta_2)u = \Delta_1 u+\Delta_2 u = Q_1[u]+Q_2[u]$ .

(*c*) $\Delta_1+\Delta_2$ commutes with total derivatives $D_i$ , since so do $\Delta_1$ and $\Delta_2$ .

(*d*) Given that each of $\Delta_1$ and $\Delta_2$ satisfies the Leibniz rule (9), it is not hard to show that the same is true for $\Delta_1+\Delta_2$ .

*Proposition 3.* Let $\Delta_1$ and $\Delta_2$ be characteristic derivatives with respect to the characteristics $Q_1[u]$ and $Q_2[u]$, respectively; i.e., $\Delta_1 u=Q_1[u]$, $\Delta_2 u=Q_2[u]$. We define the operator $[\Delta_1 , \Delta_2]$ (*Lie bracket* of $\Delta_1$ and $\Delta_2$) by

$$[\Delta_1 , \Delta_2]\,F[u] \equiv \Delta_1\,(\Delta_2 F[u]) - \Delta_2\,(\Delta_1 F[u])\,.$$

Then, $[\Delta_1 , \Delta_2]$ is a characteristic derivative with characteristic $\Delta_1 Q_2[u]-\Delta_2 Q_1[u]$. That is,

$$[\Delta_1, \Delta_2] = \Delta_Q \quad \text{with} \quad Q[u] = \Delta_1 Q_2[u] - \Delta_2 Q_1[u] \equiv Q_{1,2}[u] \tag{16}$$

*Proof.* (*a*) The linearity of $[\Delta_1 , \Delta_2]$ follows from the linearity of $\Delta_1$ and $\Delta_2$ .

(*b*) For $F[u]=u$, $[\Delta_1 , \Delta_2]u = \Delta_1\,(\Delta_2 u) - \Delta_2\,(\Delta_1 u) = \Delta_1 Q_2[u]-\Delta_2 Q_1[u] \equiv Q_{1,2}[u]$ .

(*c*) $[\Delta_1 , \Delta_2]$ commutes with total derivatives $D_i$ , since so do $\Delta_1$ and $\Delta_2$ .

(*d*) Given that each of $\Delta_1$ and $\Delta_2$ satisfies the Leibniz rule (9), one can show (after some algebra and cancellation of terms) that the same is true for $[\Delta_1 , \Delta_2]$.

In the case where $u$ (thus the $Q$'s also) is scalar-valued, the Lie bracket admits a standard differential-operator representation:

$$[\Delta_1, \Delta_2] = Q_{1,2}[u]\frac{\partial}{\partial u} + \left(D_i Q_{1,2}\right)\frac{\partial}{\partial u_i} + \left(D_i D_j Q_{1,2}\right)\frac{\partial}{\partial u_{ij}} + \cdots \qquad (17)$$

where $Q_{1,2}[u] = [\Delta_1 , \Delta_2]\, u = \Delta_1 Q_2[u] - \Delta_2 Q_1[u]$ .

The Lie bracket $[\Delta_1 , \Delta_2]$ has the following properties:

1. $[\Delta_1 , a\Delta_2 + b\Delta_3] = a[\Delta_1 , \Delta_2] + b[\Delta_1 , \Delta_3]$ ;
   $[a\Delta_1 + b\Delta_2 , \Delta_3] = a[\Delta_1 , \Delta_3] + b[\Delta_2 , \Delta_3]$ $\quad$ ($a$, $b$ = const.)
2. $[\Delta_1 , \Delta_2] = -[\Delta_2 , \Delta_1]$ $\quad$ (*antisymmetry*)
3. $[\Delta_1 , [\Delta_2 , \Delta_3]] + [\Delta_2 , [\Delta_3 , \Delta_1]] + [\Delta_3 , [\Delta_1 , \Delta_2]] = 0$ ;
   $[[\Delta_1 , \Delta_2] , \Delta_3] + [[\Delta_2 , \Delta_3] , \Delta_1] + [[\Delta_3 , \Delta_1] , \Delta_2] = 0$ $\quad$ (*Jacobi identity*)

## 4. Infinitesimal symmetry transformations of a PDE

Let $F[u]=0$ be a PDE in the independent variables $x^k \equiv x^1, x^2, \ldots$ , and the (generally) matrix-valued dependent variable $u$. A transformation $u(x^k) \rightarrow u'(x^k)$ from the function $u$ to a new function $u'$ represents a *symmetry* of the PDE if the following condition is satisfied: $u'(x^k)$ is a solution of $F[u]=0$ *when* $u(x^k)$ is a solution; that is, $F[u']=0$ when $F[u]=0$.

We will restrict our attention to *continuous symmetries* and, for the moment, to infinitesimal transformations. Although such symmetries may involve transformations of the independent variables ($x^k$), they may equivalently be expressed as transformations of $u$ alone (see [1], Chap. 5), i.e., as transformations in the fiber space.

An infinitesimal symmetry transformation is written symbolically as

$$u \rightarrow u' = u + \delta u$$

where $\delta u$ is an infinitesimal quantity, in the sense that all powers $(\delta u)^n$ with $n>1$ may be neglected. The *symmetry condition* is thus written

$$F[u+\delta u] = 0 \quad \text{when} \quad F[u] = 0 \qquad (18)$$

An infinitesimal change $\delta u$ of $u$ induces a change $\delta F[u]$ of $F[u]$, where

$$\delta F[u] = F[u+\delta u] - F[u] \quad \Leftrightarrow \quad F[u+\delta u] = F[u] + \delta F[u] \qquad (19)$$

Now, if $\delta u$ is an infinitesimal symmetry and if $u$ is a solution of $F[u]=0$, then $u+\delta u$ also is a solution; that is, $F[u+\delta u]=0$. This means that $\delta F[u]=0$ when $F[u]=0$. The symmetry condition (18) is thus written as follows:

$$\delta F[u]=0 \quad \mathrm{mod} \quad F[u] \qquad (20)$$

A finite symmetry transformation (we denote it $M$) of the PDE $F[u]=0$ produces a one-parameter family of solutions of the PDE from any given solution $u(x^k)$. We express this by writing

$$M:\ u(x^k) \rightarrow \bar{u}(x^k;\alpha) \quad \text{with} \quad \bar{u}(x^k;0)=u(x^k) \qquad (21)$$

For infinitesimal values of the parameter $\alpha$,

$$\bar{u}(x^k;\alpha) \simeq u(x^k)+\alpha Q[u] \quad \text{where} \quad Q[u]=\left.\frac{d\bar{u}}{d\alpha}\right|_{\alpha=0} \qquad (22)$$

The function $Q[u] \equiv Q(x^k, u, u_k, u_{kl}, ...)$ in the jet space is called the *characteristic* of the symmetry (or, the *symmetry characteristic*). Putting

$$\delta u=\bar{u}(x^k;\alpha)-u(x^k) \qquad (23)$$

we write, for infinitesimal $\alpha$,

$$\delta u = \alpha Q[u] \qquad (24)$$

We notice that the infinitesimal operator $\delta$ has the following properties:

1. According to its definition (19), $\delta$ is a linear operator:

$$\delta(F[u]+G[u]) = (F[u+\delta u]+ G[u+\delta u]) - (F[u]+G[u]) = \delta F[u]+\delta G[u] \ .$$

2. By the nature of our symmetry transformations, $\delta$ produces changes in the fiber space while it doesn't affect functions $f(x^k)$ in the base space [this is implicitly stated in (23)].

3. Since $\delta$ represents a difference, it commutes with all total derivatives $D_i$:

$$\delta\,(D_i A[u]) = D_i\,(\delta A[u]) \ .$$

In particular, for $A[u]=u$,

$$\delta u_i = \delta\,(D_i\, u) = D_i\,(\delta u) = \alpha D_i\, Q[u]\ ,$$

where we have used (24).

4. Since $\delta$ expresses an infinitesimal change, it may be approximated to a differential; in particular, it satisfies the Leibniz rule:

$$\delta\,(A[u]B[u]) = (\delta A[u])\,B[u] + A[u]\,\delta B[u] \ .$$

For example, $\delta(u^2) = \delta(uu) = (\delta u)u+u\delta u = \alpha\,(Qu+uQ)$ .

Now, consider the characteristic derivative $\Delta_Q$ with respect to the symmetry characteristic $Q[u]$. According to (7),

$$\Delta_Q u = Q[u] \tag{25}$$

We observe that the infinitesimal operator $\delta$ and the characteristic derivative $\Delta_Q$ share common properties. From (24) and (25) it follows that the two linear operators are related by

$$\delta u = \alpha \Delta_Q u \tag{26}$$

and, by extension,

$$\delta u_i = \alpha D_i Q[u] = \alpha \Delta_Q u_i \,, \ \text{etc.}$$

[see (10)]. Moreover, for scalar-valued $u$ and by the infinitesimal character of the operator $\delta$, we may write:

$$\delta F[u] = \frac{\partial F}{\partial u}\delta u + \frac{\partial F}{\partial u_i}\delta u_i + \cdots = \alpha\left(\frac{\partial F}{\partial u}Q[u] + \frac{\partial F}{\partial u_i}D_i Q[u] + \frac{\partial F}{\partial u_{ij}}D_i D_j Q[u] + \cdots\right)$$

while, by (13),

$$\Delta_Q F[u] = \frac{\partial F}{\partial u}Q[u] + \frac{\partial F}{\partial u_i}D_i Q[u] + \frac{\partial F}{\partial u_{ij}}D_i D_j Q[u] + \cdots \tag{27}$$

The above observations lead us to the conclusion that, in general, the following relation is true:

$$\delta F[u] = \alpha \Delta_Q F[u] \tag{28}$$

The symmetry condition (20) is thus written:

$$\Delta_Q F[u] = 0 \mod F[u] \tag{29}$$

In particular, for a scalar-valued $u$ the above condition assumes the form

$$\frac{\partial F}{\partial u}Q[u] + \frac{\partial F}{\partial u_i}D_i Q[u] + \frac{\partial F}{\partial u_{ij}}D_i D_j Q[u] + \cdots = 0 \mod F[u] \tag{30}$$

which is a linear PDE for $Q[u]$. More generally, for matrix-valued $u$ and for a function $F[u]$ of the form (1), the symmetry condition for the PDE $F[u]=0$ is a *linear* PDE for the symmetry characteristic $Q[u]$. We write this PDE symbolically as

$$S(Q\,;u) \equiv \Delta_Q F[u] = 0 \mod F[u] \tag{31}$$

where the function $S(Q\,;u)$ is linear in $Q$ and all total derivatives of $Q$. The linearity of $S$ in $Q$ follows from the Leibniz rule and the specific form (1) of $F[u]$.

Below is a list of formulas that may be useful in calculations:

- $\Delta_Q u_i = D_i Q[u]$ , $\Delta_Q u_{ij} = D_i D_j Q[u]$ , etc.
- $\Delta_Q u^2 = \Delta_Q (uu) = Q[u]u + uQ[u]$ (etc.)
- $\Delta_Q (u^{-1}) = -u^{-1} (\Delta_Q u) u^{-1} = -u^{-1} Q[u] u^{-1}$
- $\Delta_Q [A[u], B[u]] = [\Delta_Q A, B] + [A, \Delta_Q B]$

*Comment:* According to (29), $\Delta_Q F[u]$ vanishes if $F[u]$ vanishes. Given that $\Delta_Q$ is a linear operator, the reader may wonder whether this condition is identically satisfied (a linear operator acting on a zero function always produces a zero function!). Note, however, that the function $F[u]$ is *not identically* zero; it becomes zero only *for solutions* of the given PDE. What we need to do, therefore, is to first evaluate $\Delta_Q F[u]$ for *arbitrary* $u$ and *then* demand that the result vanish *when* $u$ is a solution of the PDE $F[u]$=0.

An alternative – and perhaps more transparent – version of the symmetry condition (29) is the requirement that the following relation be satisfied:

$$\Delta_Q F[u] = \hat{L} F[u] \tag{32}$$

where $\hat{L}$ is a linear operator acting on functions in the jet space (see, e.g., [1], Chap. 2 and 5, for a rigorous justification of this condition in the case of scalar-valued PDEs). For example, one may have

$$\Delta_Q F[u] = \sum_i \beta_i (x^k) D_i F[u] + \sum_{i,j} \gamma_{ij} (x^k) D_i D_j F[u] + A(x^k) F[u] + F[u] B(x^k)$$

where the $\beta_i$ and $\gamma_{ij}$ are scalar-valued while $A$ and $B$ are matrix-valued. Let us see some examples, restricting for the moment our attention to scalar PDEs.

*Example 1.* The *sine-Gordon (s-G) equation* is written

$$F[u] \equiv u_{xt} - \sin u = 0 \; .$$

Here, $(x^1, x^2) \equiv (x, t)$. Since $\sin u$ can be expanded into an infinite series in powers of $u$, we see that $F[u]$ has the required form (1). Moreover, since $u$ is a scalar function, we can write the symmetry condition by using (30):

$$S(Q; u) \equiv Q_{xt} - (\cos u)\, Q = 0 \quad \text{mod} \quad F[u]$$

where $S(Q; u) = \Delta_Q F[u]$ and where by subscripts we denote total differentiations with respect to the indicated variables. Let us verify the solution $Q[u] = u_x$ . As will be shown in Sec. 9, this characteristic produces the finite symmetry transformation

$$M: \; u(x,t) \rightarrow \bar{u}(x,t;\alpha) = u(x+\alpha, t) \tag{33}$$

which implies that, if $u(x,t)$ is a solution of the s-G equation, then $\bar{u}(x,t) = u(x+\alpha, t)$ also is a solution. We have:

$$Q_{xt}-(\cos u)\,Q=(u_x)_{xt}-(\cos u)\,u_x=(u_{xt}-\sin u)_x=D_x F[u]=0 \mod F[u]\ .$$

Notice that $\Delta_Q F[u]$ is of the form (32), with $\hat{L}\equiv D_x$. Similarly, the characteristic $Q[u]=u_t$ corresponds to the symmetry

$$M:\ u(x,t)\rightarrow \bar{u}(x,t;\alpha)=u(x,\,t+\alpha) \tag{34}$$

That is, if $u(x,t)$ is a solution of the s-G equation, then so is $\bar{u}(x,t)=u(x,\,t+\alpha)$. The symmetries (33) and (34) reflect the fact that the s-G equation does not contain the variables $x$ and $t$ explicitly. (Of course, this equation has many more symmetries that are not displayed here; see, e.g., [1].)

*Example 2.* The *heat equation* is written

$$F\,[u]\equiv u_t-u_{xx}=0\ .$$

The symmetry condition (30) reads

$$S(Q\,;u)\equiv Q_t-Q_{xx}=0 \mod F[u]$$

where $S(Q;u)=\Delta_Q F[u]$. As is easy to show, the symmetries (33) and (34) are valid here, too. Let us now try the solution $Q[u]=u$. We have:

$$Q_t-Q_{xx}=\ u_t-u_{xx}=\ F[u]=0 \mod F[u]\ .$$

As will be shown in Sec. 9, this symmetry corresponds to the transformation

$$M:\ u(x,t)\rightarrow \bar{u}(x,t;\alpha)=e^{\alpha}\,u(x,t) \tag{35}$$

and is a consequence of the linearity of the heat equation.

*Example 3.* One form of the *Burgers equation* is

$$F[u]\equiv u_t-u_{xx}-u_x^{\,2}=0\ .$$

The symmetry condition (30) is written

$$S(Q\,;u)\equiv Q_t-Q_{xx}-2u_xQ_x=0 \mod F[u]$$

where, as always, $S(Q;u)=\Delta_Q F[u]$. Putting $Q=u_x$ and $Q=u_t$, we verify the symmetries (33) and (34):

$$Q_t-Q_{xx}-2u_xQ_x=u_{xt}-u_{xxx}-2u_xu_{xx}=D_xF[u]=0 \mod F[u]$$
$$Q_t-Q_{xx}-2u_xQ_x=u_{tt}-u_{xxt}-2u_xu_{xt}=D_tF[u]=0 \mod F[u]$$

Note again that $\Delta_Q F[u]$ is of the form (32), with $\hat{L}\equiv D_x$ and $\hat{L}\equiv D_t$. Another symmetry is $Q\,[u]=1$, which corresponds to the transformation (see Sec. 9)

$$M:\; u(x,t) \to \bar{u}(x,t;\alpha) = u(x,t) + \alpha \qquad (36)$$

and is a consequence of the fact that $u$ enters $F[u]$ only through its derivatives.

*Example 4.* The *wave equation* is written

$$F[u] \equiv u_{tt} - c^2 u_{xx} = 0 \quad (c = \text{const.})$$

and its symmetry condition reads

$$S(Q;u) \equiv Q_{tt} - c^2 Q_{xx} = 0 \mod F[u]\,.$$

The solution $Q[u] = x\,u_x + t\,u_t$ corresponds to the symmetry transformation (Sec. 9)

$$M:\; u(x,t) \to \bar{u}(x,t;\alpha) = u(e^{\alpha}x,\, e^{\alpha}t) \qquad (37)$$

expressing the invariance of the wave equation under a scale change of $x$ and $t$. [The reader may show that the transformations (33) – (36) also express symmetries of the wave equation.]

It is remarkable that each of the above PDEs admits an infinite set of symmetry transformations [1]. An effective method for finding such infinite sets is the use of a *recursion operator*, which produces a new symmetry characteristic every time it acts on a known characteristic. More on recursion operators will be said in Sec. 6.

## 5. The Lie algebra of symmetries

As is well known [1] the set of symmetries of a PDE $F[u]=0$ has the structure of a Lie algebra. Let us demonstrate this property in the context of our formalism.

*Proposition.* Let $\mathcal{L}$ be the set of characteristic derivatives $\Delta_Q$ with respect to the symmetry characteristics $Q[u]$ of the PDE $F[u]=0$. The set $\mathcal{L}$ is a (finite or infinite-dimensional) Lie subalgebra of the Lie algebra of characteristic derivatives acting on functions $F[u]$ in the jet space (cf. Sec. 3).

*Proof.* (*a*) Let $\Delta_Q \in \mathcal{L} \Leftrightarrow \Delta_Q F[u]=0 \pmod{F[u]}$. If $\lambda$ is a constant then $(\lambda\Delta_Q)F[u] \equiv \lambda\Delta_Q F[u]=0$, which means that $\lambda\Delta_Q \in \mathcal{L}$. Given that $\lambda\Delta_Q = \Delta_{\lambda Q}$ [see Eq. (14)] we conclude that, if $Q[u]$ is a symmetry characteristic of $F[u]=0$, then so is $\lambda Q[u]$.

(*b*) Let $\Delta_1 \in \mathcal{L}$ and $\Delta_2 \in \mathcal{L}$ be characteristic derivatives with respect to the symmetry characteristics $Q_1[u]$ and $Q_2[u]$, respectively. Then, $\Delta_1 F[u]=0$, $\Delta_2 F[u]=0$, and hence $(\Delta_1+\Delta_2)F[u] \equiv \Delta_1 F[u]+\Delta_2 F[u]=0$; therefore, $(\Delta_1+\Delta_2) \in \mathcal{L}$. It also follows from Eq. (15) that, if $Q_1[u]$ and $Q_2[u]$ are symmetry characteristics of $F[u]=0$, then so is their sum $Q_1[u]+Q_2[u]$.

(*c*) Let $\Delta_1 \in \mathcal{L}$ and $\Delta_2 \in \mathcal{L}$, as above. Then, by (32),

$$\Delta_1 F[u] = \hat{L}_1 F[u]\,, \quad \Delta_2 F[u] = \hat{L}_2 F[u]\,.$$

Now, by the definition of the Lie bracket and by the linearity of both $\Delta_i$ and $\hat{L}_i$ ($i$=1,2) we have:

$$[\Delta_1, \Delta_2]F[u] = \Delta_1(\Delta_2 F[u]) - \Delta_2(\Delta_1 F[u]) = \Delta_1(\hat{L}_2 F[u]) - \Delta_2(\hat{L}_1 F[u])$$
$$\equiv (\Delta_1 \hat{L}_2 - \Delta_2 \hat{L}_1)F[u] = 0 \mod F[u]$$

We thus conclude that $[\Delta_1, \Delta_2] \in \mathcal{L}$. Moreover, it follows from Eq. (16) that, if $Q_1[u]$ and $Q_2[u]$ are symmetry characteristics of $F[u]=0$, then so is the function

$$Q_{1,2}[u] = \Delta_1 Q_2[u] - \Delta_2 Q_1[u] .$$

Assume now that the PDE $F[u]=0$ has an $n$-dimensional symmetry algebra $\mathcal{L}$ (which may be a finite subalgebra of an infinite-dimensional symmetry Lie algebra). Let $\{\Delta_1, \Delta_2, \ldots, \Delta_n\} \equiv \{\Delta_k\}$, with corresponding symmetry characteristics $\{Q_k\}$, be a set of $n$ linearly independent operators that constitute a basis of $\mathcal{L}$, and let $\Delta_i$, $\Delta_j$ be any two elements of this basis. Given that $[\Delta_i, \Delta_j] \in \mathcal{L}$, this Lie bracket must be expressible as a linear combination of the $\{\Delta_k\}$, with constant coefficients. We write

$$[\Delta_i, \Delta_j] = \sum_{k=1}^{n} c_{ij}^k \Delta_k \tag{38}$$

where the coefficients of the $\Delta_k$ in the sum are the antisymmetric *structure constants* of the Lie algebra $\mathcal{L}$ in the basis $\{\Delta_k\}$.

The operator relation (38) can be expressed in an equivalent, characteristic form by allowing the operators on both sides to act on $u$ and by using the fact that $\Delta_k u = Q_k[u]$:

$$[\Delta_i, \Delta_j]u = \left( \sum_{k=1}^{n} c_{ij}^k \Delta_k \right) u = \sum_{k=1}^{n} c_{ij}^k (\Delta_k u) \Rightarrow$$

$$\Delta_i Q_j[u] - \Delta_j Q_i[u] = \sum_{k=1}^{n} c_{ij}^k Q_k[u] \tag{39}$$

*Example.* One of the several forms of the *Korteweg-de Vries (KdV) equation* is

$$F[u] \equiv u_t + u u_x + u_{xxx} = 0 .$$

The symmetry condition (31) is written

$$S(Q;u) \equiv Q_t + Q u_x + u Q_x + Q_{xxx} = 0 \mod F[u] \tag{40}$$

where $S(Q;u) = \Delta_Q F[u]$. The KdV equation admits a symmetry Lie algebra of infinite dimensions [1]. This algebra has a finite, 4-dimensional subalgebra $\mathcal{L}$ of *point transformations*. A symmetry operator (characteristic derivative) $\Delta_Q$ is determined by its corresponding characteristic $Q[u] = \Delta_Q u$. Thus, a basis $\{\Delta_1, \ldots, \Delta_4\}$ of $\mathcal{L}$ corresponds to a set of four independent characteristics $\{Q_1 ,..., Q_4\}$. Such a basis of characteristics is the following:

$$Q_1[u]= u_x\,,\quad Q_2[u]= u_t\,,\quad Q_3[u]= tu_x-1\,,\quad Q_4[u]= xu_x +3tu_t + 2u$$

The $Q_1 ,..., Q_4$ satisfy the PDE (40), since, as we can show,

$$S(Q_1 ; u) = D_x F[u]\,,\quad S(Q_2 ; u) = D_t F[u]\,,\quad S(Q_3 ; u) = tD_x F[u]\,,$$
$$S(Q_4 ; u) = (5 + xD_x + 3tD_t)F[u]$$

[Note once more that $\Delta_Q F[u]$ is of the form (32) in each case.] Let us now see two examples of calculating the structure constants of $\mathcal{L}$ by application of (39). We have:

$$\Delta_1 Q_2 - \Delta_2 Q_1 = \Delta_1 u_t - \Delta_2 u_x = (\Delta_1 u)_t - (\Delta_2 u)_x = (Q_1)_t - (Q_2)_x = (u_x)_t - (u_t)_x = 0$$
$$\equiv \sum_{k=1}^{4} c_{12}^k Q_k$$

Since the $Q_k$ are linearly independent, we must necessarily have $c_{12}^k = 0$, $k = 1,2,3,4$. Also,

$$\Delta_2 Q_3 - \Delta_3 Q_2 = \Delta_2 (t u_x - 1) - \Delta_3 u_t = t(\Delta_2 u)_x - (\Delta_3 u)_t = t(Q_2)_x - (Q_3)_t$$
$$= t u_{tx} - (u_x + t u_{xt}) = -u_x = -Q_1 \equiv \sum_{k=1}^{4} c_{23}^k Q_k$$

Therefore, $c_{23}^1 = -1$, $c_{23}^2 = c_{23}^3 = c_{23}^4 = 0$.

## 6. Recursion operators

Let $\delta u = \alpha Q[u]$ be an infinitesimal symmetry of the PDE $F[u]=0$, where $Q[u]$ is the symmetry characteristic. For any solution $u(x^k)$ of this PDE the function $Q[u]$ satisfies the linear PDE

$$S(Q ; u) \equiv \Delta_Q F[u] = 0 \qquad (41)$$

Because of the linearity of (41) in $Q$, the sum $Q_1[u]+Q_2[u]$ of two solutions of this PDE, as well as the multiple $\lambda Q[u]$ of any solution by a constant, also are solutions of (41) for a given $u$. Thus, for any solution $u$ of $F[u]=0$, the solutions $\{Q[u]\}$ of the PDE (41) form a linear space, which we call $S_u$.

A *recursion operator* $\hat{R}$ is a linear operator that maps the space $S_u$ into itself. Thus, if $Q[u]$ is a symmetry characteristic of $F[u]=0$ [i.e., a solution of (41) for a given $u$] then so is $\hat{R}Q[u]$:

$$S(\hat{R}Q ; u) = 0 \quad \text{when} \quad S(Q ; u) = 0 \qquad (42)$$

Obviously, any power of a recursion operator also is a recursion operator. Thus, starting with any symmetry characteristic $Q[u]$, one may in principle obtain an infinite set of such characteristics by repeated application of the recursion operator.

A new approach to recursion operators was suggested in the early 1990s [11,15-17] (see also [8-10]) according to which a recursion operator may be viewed as an *auto-Bäcklund transformation* (BT) [18,19] for the symmetry condition (41) of the PDE $F[u]=0$. By integrating the BT, one obtains new solutions $Q'[u]$ of the linear PDE (41) from known ones, $Q[u]$. Typically, this type of recursion operator produces *nonlocal* symmetries in which the symmetry characteristic depends on *integrals* (rather than derivatives) of $u$. This approach proved to be particularly effective for matrix-valued PDEs such as the nonlinear self-dual Yang-Mills equation, of which new infinite-dimensional sets of "potential symmetries" were discovered [9,11,15]. More on this will be said in Sec. 10.

## 7. An example: The chiral field equation

Let us consider the *chiral field equation*

$$F[g] \equiv (g^{-1}g_x)_x + (g^{-1}g_t)_t = 0 \tag{43}$$

where, in general, subscripts $x$ and $t$ denote total derivatives $D_x$ and $D_t$ , respectively, and where $g$ is a $GL(n,C)$-valued function of $x$ and $t$, i.e., a complex, non-singular $(n\times n)$ matrix function, differentiable for all $x$ and $t$. Let $\delta g=\alpha Q[g]$ be an infinitesimal symmetry transformation for the PDE (43), with symmetry characteristic $Q[g]=\Delta_Q\, g$. It is convenient to put

$$Q[g] = g\,\Phi[g] \quad \Leftrightarrow \quad \Phi[g] = g^{-1}Q[g]\ .$$

The symmetry condition for (43) is

$$\Delta_Q F[g] = 0 \mod F[g]\ .$$

This condition will yield a linear PDE for $Q$ or, equivalently, a linear PDE for $\Phi$. By using the properties of the characteristic derivative we find the latter PDE to be

$$S(\Phi;g) \equiv D_x\left(\Phi_x + [g^{-1}g_x, \Phi]\right) + D_t\left(\Phi_t + [g^{-1}g_t, \Phi]\right) = 0 \mod F[g] \tag{44}$$

where the square brackets denote commutators of matrices. Note that, since (44) is to be valid on solutions $g$ of (43), it can equivalently be expressed as

$$S(\Phi;g) \equiv \Phi_{xx} + \Phi_{tt} + [g^{-1}g_x\, ,\, \Phi_x] + [g^{-1}g_t\, ,\, \Phi_t] = 0 \mod F[g]\ .$$

A useful identity that will be needed in the sequel is the following:

$$(g^{-1}g_t)_x - (g^{-1}g_x)_t + [g^{-1}g_x\, ,\, g^{-1}g_t] = 0 \tag{45}$$

Let us first consider symmetry transformations in the base space, i.e., coordinate transformations of $x$, $t$. An obvious symmetry is $x$-translation, $x'=x+\alpha$, given that the PDE (43) does not contain the independent variable $x$ explicitly. For infinitesimal values of the parameter $\alpha$, we write $\delta x=\alpha$. The symmetry characteristic is $Q[g]=g_x$, so that $\Phi[g]=g^{-1}g_x$. By substituting this expression for $\Phi$ into the symmetry condition (44) and by using the identity (45), we can verify that (44) is indeed satisfied:

$$S(\Phi;g)=D_x F[g]=0 \mod F[g] \,.$$

Similarly, for $t$-translation, $t'=t+\alpha$ (infinitesimally, $\delta t=\alpha$) with $Q[g]=g_t$, we find

$$S(\Phi;g)=D_t F[g]=0 \mod F[g] \,.$$

Another obvious symmetry of (43) is a scale change of both $x$ and $t$: $x'=\lambda x$, $t'=\lambda t$. Setting $\lambda=1+\alpha$, where $\alpha$ is infinitesimal, we write $\delta x=\alpha x$, $\delta t=\alpha t$. The symmetry characteristic is $Q[g]=xg_x+tg_t$, so that $\Phi[g]=xg^{-1}g_x+tg^{-1}g_t$. Substituting for $\Phi$ into the symmetry condition (44) and using the identity (45) where necessary, we find that

$$S(\Phi;g)=(2+xD_x+tD_t)\,F[g]=0 \mod F[g] \,.$$

Let us call $Q_1[g]=g_x$, $Q_2[g]=g_t$, $Q_3[g]=xg_x+tg_t$, and let us consider the corresponding characteristic derivative operators $\Delta_i$ defined by $\Delta_i\, g=Q_i$ ($i$=1,2,3). It is then straightforward to verify the following commutation relations:

$$[\Delta_1,\Delta_2]\,g=\Delta_1 Q_2-\Delta_2 Q_1=0 \;\Leftrightarrow\; [\Delta_1,\Delta_2]=0\,;$$

$$[\Delta_1,\Delta_3]\,g=\Delta_1 Q_3-\Delta_3 Q_1=-g_x=-Q_1=-\Delta_1 g \;\Leftrightarrow\; [\Delta_1,\Delta_3]=-\Delta_1\,;$$

$$[\Delta_2,\Delta_3]\,g=\Delta_2 Q_3-\Delta_3 Q_2=-g_t=-Q_2=-\Delta_2 g \;\Leftrightarrow\; [\Delta_2,\Delta_3]=-\Delta_2\,.$$

Next, we consider the "internal" transformation (i.e., transformation in the fiber space) $g'=g\Lambda$, where $\Lambda$ is a non-singular constant matrix. Then,

$$F[g']=\Lambda^{-1}F[g]\,\Lambda=0 \mod F[g] \,,$$

which indicates that this transformation is a symmetry of (43). Setting $\Lambda=1+\alpha M$, where $\alpha$ is an infinitesimal parameter while $M$ is a constant matrix, we write, in infinitesimal form, $\delta g=\alpha gM$. The symmetry characteristic is $Q[g]=gM$, so that $\Phi[g]=M$. Substituting for $\Phi$ into the symmetry condition (44) we find:

$$S(\Phi;g)=\big[F[g]\,,M\big]=0 \mod F[g] \,.$$

Given a matrix function $g(x,t)$ satisfying the PDE (43), consider the following system of PDEs for two functions $\Phi[g]$ and $\Phi'[g]$:

$$\begin{aligned} \Phi'_x &= \Phi_t + [g^{-1}g_t\,,\,\Phi] \\ -\Phi'_t &= \Phi_x + [g^{-1}g_x\,,\,\Phi] \end{aligned} \qquad (46)$$

The integrability condition $(\Phi'_x)_t = (\Phi'_t)_x$ of this system requires that $\Phi$ satisfy the symmetry condition (44); i.e., $S(\Phi\,;\,g)=0$. Conversely, by applying the integrability condition $(\Phi_t)_x = (\Phi_x)_t$ and by using the fact that $g$ is a solution of $F[g]=0$, one finds that $\Phi'$ must also satisfy (44); i.e., $S(\Phi';\,g) = 0$.

We conclude that, for any function $g(x,t)$ satisfying the PDE (43), the system (46) is an *auto-Bäcklund transformation* (BT) [18,19] relating solutions $\Phi$ and $\Phi'$ of the symmetry condition (44) of this PDE; that is, relating different symmetries of the chiral field equation. Thus, if a symmetry characteristic $Q=g\Phi$ of the PDE (43) is known, a new characteristic $Q'=g\Phi'$ may be found by integrating the BT (46); the converse is also true. Since the BT (46) produces new symmetries from old ones it may be regarded as a *recursion operator* for the PDE (43) [8-11,15-17].

As an example, consider the internal-symmetry characteristic $Q[g]=gM$ (where $M$ is a constant matrix) corresponding to $\Phi[g]=M$. By integrating the BT (46) for $\Phi'$, we get $\Phi'=[X, M]$ and thus $Q'=g[X, M]$, where $X$ is the "potential" of the PDE (43), defined by the system of PDEs

$$X_x = g^{-1} g_t \;\; , \quad -X_t = g^{-1} g_x \tag{47}$$

Note the *nonlocal* character of the BT-produced symmetry $Q'$, due to the presence of the potential $X$. Indeed, as seen from (47), in order to find $X$ one has to *integrate* the chiral field $g$ with respect to the independent variables $x$ and $t$. The above process can be continued indefinitely by repeated application of the recursion operator (46), leading to an infinite sequence of increasingly nonlocal symmetries.

Unfortunately, as the reader may check, no new information is furnished by the BT (46) in the case of coordinate symmetries (for example, by applying the BT for $Q=g_x$ we get the known symmetry $Q'=g_t$ ). A recursion operator of the form (46), however, does produce new nonlocal symmetries from coordinate symmetries in related problems with more than two independent variables, such as the self-dual Yang-Mills equation [8-11,15].

## 8. Generation of finite symmetry transformations

As we saw in Sec. 4, given a symmetry operator $\Delta_Q$ one may immediately define an infinitesimal symmetry of a PDE. Our starting point, however, was the idea of using a *finite* symmetry transformation to generate a one-parameter family of solutions of the PDE. We thus need to generalize the discussion of Sec. 4 by allowing the parameter $\alpha$ to assume finite values.

According to (7), the characteristic derivative $\Delta_Q$ with respect to the characteristic function $Q[u]$ satisfies the relation

$$\Delta_{Q[u]}\, u = Q[u] \tag{48}$$

By (48) and the properties of $\Delta_Q$ one may determine the action of $\Delta_Q$ on any function $F[u]$ of the form (1), thus construct a new function $\Delta_{Q[u]}F[u]$.

Obviously, a change of $u$ will induce a corresponding change on any function $F[u]$. Given a function $u(x^k)$, a continuous change of $u$ may be expressed as a one-parameter family of transformations

$$M:\ u(x^k) \to \bar{u}(x^k;\alpha) \quad \text{with} \quad \bar{u}(x^k;0) = u(x^k) \tag{49}$$

where $\alpha \geq 0$. We suppress the independent variables $x^k$, which are unaffected by the transformation $M$, and write, simply,

$$M:\ u \to \bar{u}(\alpha) \quad \text{with} \quad \bar{u}(0) = u\,.$$

Expanding $\bar{u}(\alpha)$ in powers of $\alpha$, we have:

$$\bar{u}(\alpha) = u + \alpha Q[u] + \text{higher-order terms in } \alpha \tag{50}$$

where $Q[u]$ is given by

$$\frac{d}{d\alpha}\bar{u}(\alpha)\Big|_{\alpha=0} = Q[u] = \Delta_{Q[u]}u \tag{51}$$

Now, we assume that, for finite values of the parameter $\alpha$,

$$\frac{d}{d\alpha}\bar{u}(\alpha) = Q\big[\bar{u}(\alpha)\big] = \Delta_{Q[\bar{u}(\alpha)]}\,\bar{u}(\alpha) \tag{52}$$

which is obviously consistent with (51). By the properties of the characteristic derivative it then follows that, for any function $F[u]$ of the form (1),

$$\frac{d}{d\alpha}F\big[\bar{u}(\alpha)\big] = \Delta_{Q[\bar{u}(\alpha)]}\,F\big[\bar{u}(\alpha)\big] \tag{53}$$

As an example, let $F[u] = u^2 \Rightarrow F[\bar{u}(\alpha)] = \bar{u}(\alpha)\bar{u}(\alpha)$. Then, by (52) and by using the Leibniz rule, we have:

$$\begin{aligned}\frac{d}{d\alpha}F[\bar{u}(\alpha)] &= \frac{d\bar{u}(\alpha)}{d\alpha}\,\bar{u}(\alpha) + \bar{u}(\alpha)\frac{d\bar{u}(\alpha)}{d\alpha}\\ &= \left\{\Delta_{Q[\bar{u}(\alpha)]}\bar{u}(\alpha)\right\}\bar{u}(\alpha) + \bar{u}(\alpha)\Delta_{Q[\bar{u}(\alpha)]}\bar{u}(\alpha)\\ &= \Delta_{Q[\bar{u}(\alpha)]}\big(\bar{u}(\alpha)\bar{u}(\alpha)\big) = \Delta_{Q[\bar{u}(\alpha)]}\,F\big[\bar{u}(\alpha)\big]\end{aligned}$$

Equation (52), together with the initial condition $\bar{u}(0) = u$, is an initial-value problem that, upon integration, yields a one-parameter transformation of the form (49). We say that the operator $\Delta_Q$ is the *generator* of this transformation. As regards its action on functions, the operator $\Delta_Q$ is seen to be equivalent to the *Lie derivative* of differential geometry (see, e.g., Chap. 5 of [19]). And, the latter derivative plays a key role in the differential-geometric methods for studying invariance properties of PDEs [2-4]. We now revisit the symmetry problem for PDEs in the context of our present, more abstract algebraic formalism.

The transformation $M$ in Eq. (49) is a *symmetry transformation* for the PDE $F[u]=0$ if it leaves this PDE invariant, in the sense that $F[\bar{u}(\alpha)] = 0$ *if* $F[u] = 0$. We write:

$$F[\bar{u}(\alpha)]=0 \mod F[u] \tag{54}$$

So, if $\bar{u}(x^k;0)=u(x^k)$ is a solution of $F[u]=0$, then so is $\bar{u}(x^k;\alpha)$ *for all values of the parameter* $\alpha>0$. This means that $F[\bar{u}(\alpha)]$, viewed as a function of $\alpha$, retains a constant (zero) value under continuous changes of $\alpha$. In mathematical terms,

$$\frac{d}{d\alpha}F\left[\bar{u}(\alpha)\right]=0 \mod F\left[\bar{u}(\alpha)\right]$$

or, in view of (53),

$$\Delta_{Q[\bar{u}(\alpha)]}F\left[\bar{u}(\alpha)\right]=0 \mod F\left[\bar{u}(\alpha)\right] .$$

Since this must be valid for any value of $\alpha$, the above relation will still be true if we replace $\alpha$ by a new parameter $\beta=\alpha+c$, where $c$ is any constant such that $\beta\geq0$. In particular, by choosing $c=-\alpha \Rightarrow \beta=0$, we rewrite the above equation in the simpler form

$$\Delta_{Q[u]}F[u]=0 \mod F[u] \tag{55}$$

which is the condition for invariance of the PDE $F[u]=0$. As we have seen, this condition yields a linear PDE for the symmetry characteristic $Q[u]$, of the form

$$S(Q;u)\equiv\Delta_{Q[u]}F[u]=0 \mod F[u] \tag{56}$$

where the expression $S(Q;u)$ is linear in $Q$ and all total derivatives of $Q$. In particular, for scalar-valued $u$ (thus scalar $Q[u]$ also) the operator $\Delta_Q$ has the form (13) and the symmetry condition (56) takes on the form

$$\frac{\partial F}{\partial u}Q[u]+\frac{\partial F}{\partial u_i}D_iQ[u]+\frac{\partial F}{\partial u_{ij}}D_iD_jQ[u]+\cdots=0 \mod F[u] \tag{57}$$

An important class of symmetries is *local* (point) symmetries. As discussed in [1], the symmetry characteristic $Q[u]$ of a local symmetry cannot depend on second or higher-order derivatives of $u$ with respect to the $x^k$, while the dependence of $Q$ on the first-order derivatives $u_k$ is also subject to certain restrictions. Once a local symmetry characteristic $Q[u]$ is found by solving (56) or (57), a one-parameter family of symmetry transformations of the PDE $F[u]=0$, of the form

$$M:\ u(x^k)\rightarrow\bar{u}(x^k;\alpha)\ \ ;\ \ \ \bar{u}(x^k;0)=u(x^k) \tag{58}$$

is obtained by solving the initial-value problem [cf. Eq. (52)]

$$\begin{gathered}\frac{d}{d\alpha}\bar{u}(x^k;\alpha)=Q\left[\bar{u}\right]\\ \bar{u}(x^k;0)=u(x^k)\end{gathered} \tag{59}$$

## 9. Example: The two-dimensional Laplace equation

As an example for illustrating the process of finding one-parameter symmetry transformations of a PDE, we choose the two-dimensional *Laplace equation* for a scalar function $u(x,t)$:

$$F[u] \equiv u_{xx} + u_{tt} = 0 \tag{60}$$

Here, $(x^1, x^2) \equiv (x, t)$. The symmetry condition (56) or (57) yields the linear PDE

$$S(Q;u) \equiv Q_{xx} + Q_{tt} = 0 \quad \text{mod} \quad F[u] \tag{61}$$

where subscripts indicate total differentiations. Each symmetry of the PDE (60) corresponds to a solution $Q[u]$ of (61) and leads to a one-parameter family of symmetry transformations (58) by solving the initial-value problem (59). Let us see some examples:

1. $Q[u]=1$ is a solution of (61), hence a symmetry characteristic of (60). The initial-value problem (59) is written

$$\frac{d}{d\alpha}\bar{u}(\alpha) = 1 \;\; ; \quad \bar{u}(0) = u$$

which is easily integrated to give $\bar{u}(x,t;\alpha) = u(x,t) + \alpha$ . Thus, if $u(x,t)$ is a solution of (60), then so is $u(x,t)+\alpha$ . This symmetry reflects the fact that $u$ enters the PDE (60) only through its derivatives (i.e., $F[u]$ does not contain $u$ itself).

2. For the symmetry characteristic $Q[u]=u$ we have

$$\frac{d}{d\alpha}\bar{u}(\alpha) = \bar{u} \;\; ; \quad \bar{u}(0) = u$$

with solution $\bar{u}(x,t;\alpha) = e^{\alpha} u(x,t)$. Thus, if $u(x,t)$ is a solution of (60), then so is $\lambda u(x,t)$ for any constant $\lambda$. This symmetry reflects the fact that the PDE (60) is homogeneous linear in $u$.

3. $Q[u]=u_x$ is another symmetry characteristic; indeed, note that $S(Q;u)=D_x F[u]=0$ when $F[u]=0$. The initial-value problem is written

$$\frac{d}{d\alpha}\bar{u}(x,t;\alpha) = \bar{u}_x \;\; ; \quad \bar{u}(x,t;0) = u(x,t)\,.$$

A way to solve this is to consider the parameter $\alpha$ as a variable of equal footing with $x$ and $t$. The above ordinary differential equation then becomes a homogeneous linear first-order PDE that can be integrated by standard methods (see, e.g., Chap. 4 of [19]):

$$\bar{u}_x - \bar{u}_\alpha = 0 \tag{62}$$

We form the characteristic system

$$\frac{dx}{1} = \frac{dt}{0} = \frac{d\alpha}{-1} = \frac{d\bar{u}}{0}$$

and we seek 3 first integrals of this system. These are $\bar{u}=C_1$, $t=C_2$, $x+\alpha=C_3$. The general solution of (62) is then $\Phi(C_1, C_2, C_3)=0$, where the function $\Phi$ is arbitrary. That is,

$$\Phi(\bar{u},\, t,\, x+\alpha)=0 \;\Rightarrow\; \bar{u}(x,t;\alpha)=\omega(x+\alpha,\, t)\,.$$

By the initial condition $\bar{u}(x,t\,;0)=u(x,t)$ we have that $\omega(x,t)=u(x,t)$, and hence $\omega(x+\alpha, t)=u(x+\alpha, t)$. Thus, finally, $\bar{u}(x,t;\alpha)=u(x+\alpha,\, t)$.

In a similar manner, from the symmetry characteristic $Q[u]=u_t$ we get the transformation $\bar{u}(x,t;\alpha)=u(x,\, t+\alpha)$. The two symmetries found above reflect the fact that the PDE (60) does not contain the independent variables $x$ and $t$ explicitly.

4. For the symmetry characteristic $Q[u]=xu_x+tu_t$ we have

$$\frac{d}{d\alpha}\bar{u}(x,t\,;\alpha)=x\bar{u}_x+t\bar{u}_t \;\;;\;\; \bar{u}(x,t\,;0)=u(x,t)\,.$$

Working as in the previous example, we form the first-order PDE

$$x\bar{u}_x+t\bar{u}_t-\bar{u}_\alpha=0 \tag{63}$$

with characteristic system

$$\frac{dx}{x}=\frac{dt}{t}=\frac{d\alpha}{-1}=\frac{d\bar{u}}{0}\;.$$

Three first integrals are $\bar{u}=C_1$, $\ln x+\alpha=C_2$, $\ln t+\alpha=C_3$. The general solution of the PDE (63) is $\Phi(C_1,C_2,C_3)=0$, with $\Phi$ arbitrary. That is,

$$\Phi(\bar{u},\, \ln x+\alpha,\, \ln t+\alpha)=0 \;\Rightarrow$$
$$\bar{u}(x,t;\alpha)=\omega(\ln x+\alpha,\, \ln t+\alpha)=\omega\left[\ln\left(e^{\alpha}x\right),\, \ln\left(e^{\alpha}t\right)\right]$$

By the initial condition $\bar{u}(x,t\,;0)=u(x,t)$ we have that $\omega(\ln x, \ln t)=u(x,t)$, and hence $\omega[\ln(e^{\alpha}x), \ln(e^{\alpha}t)]=u(e^{\alpha}x, e^{\alpha}t)$. Thus, finally, $\bar{u}(x,t;\alpha)=u(e^{\alpha}x,\, e^{\alpha}t)$. This transformation expresses the invariance of the PDE (60) under a scale change $x\to\lambda x$, $t\to\lambda t$ of $x$ and $t$.

5. $Q[u]= tu_x - xu_t$ is a symmetry characteristic since $S(Q;u)=(tD_x - xD_t)F[u]=0$ when $F[u]=0$. We write

$$\frac{d}{d\alpha}\bar{u}(x,t\,;\alpha)=t\bar{u}_x-x\bar{u}_t \;\;;\;\; \bar{u}(x,t\,;0)=u(x,t)$$

and form the PDE

$$t\bar{u}_x-x\bar{u}_t-\bar{u}_\alpha=0 \tag{64}$$

with characteristic system

$$\frac{dx}{t}=\frac{dt}{-x}=\frac{d\alpha}{-1}=\frac{d\bar{u}}{0} \tag{65}$$

One first integral is $\bar{u}=C_1$. Another one is found from $dx/t = - dt/x \Rightarrow d(x^2+t^2)=0 \Rightarrow$ $x^2+t^2 = C_2$. A third integral is

$$\alpha + \arctan(x/t) = C_3\,.$$

Let us prove this. Setting $x/t=z$, we have:

$$d(\alpha+\arctan z)=d\alpha+\frac{dz}{1+z^2}=d\alpha+\frac{tdx-xdt}{x^2+t^2}.$$

But, by the system (65), $dx=-td\alpha$ and $dt=xd\alpha$, so that

$$d\left[\alpha+\arctan(x/t)\right]=d\alpha-\frac{(x^2+t^2)d\alpha}{x^2+t^2}=0, \quad \text{q.e.d.}$$

The general solution of the PDE (64) is $\Phi(C_1,C_2,C_3)=0$, with $\Phi$ arbitrary. That is,

$$\Phi\left[\bar{u},\ x^2+t^2,\ \alpha+\arctan(x/t)\right]=0 \ \Rightarrow$$

$$\bar{u}(x,t;\alpha)=\omega\left[x^2+t^2,\ \alpha+\arctan(x/t)\right] \tag{66}$$

By the initial condition $\bar{u}(x,t;0)=u(x,t)$ we have that $\omega[x^2+t^2,\ \arctan(x/t)]=u(x,t)$. Putting $x\cos\alpha+t\sin\alpha$ and $t\cos\alpha-x\sin\alpha$ in place of $x$ and $t$, respectively, we find:

$$u(x\cos\alpha+t\sin\alpha,\ t\cos\alpha-x\sin\alpha)=\omega\left(x^2+t^2,\ \arctan\frac{x\cos\alpha+t\sin\alpha}{t\cos\alpha-x\sin\alpha}\right) \tag{67}$$

On the other hand, we can show that

$$\tan\left[\alpha+\arctan(x/t)\right]=\frac{x\cos\alpha+t\sin\alpha}{t\cos\alpha-x\sin\alpha} \ \Rightarrow$$

$$\alpha+\arctan(x/t)=\arctan\frac{x\cos\alpha+t\sin\alpha}{t\cos\alpha-x\sin\alpha} \tag{68}$$

From (67) and (68) we have that

$$u(x\cos\alpha+t\sin\alpha,\ t\cos\alpha-x\sin\alpha)=\omega\,[x^2+t^2,\ \alpha+\arctan(x/t)]\,.$$

Thus (66) assumes the final form

$$\bar{u}(x,t;\alpha)=u(x\cos\alpha+t\sin\alpha,\ -x\sin\alpha+t\cos\alpha) \tag{69}$$

The transformation (69) admits a certain geometrical interpretation that becomes evident if we define the new variables $x'=x\cos\alpha+t\sin\alpha$ and $t'=-x\sin\alpha+t\cos\alpha$. In matrix form,

$$\begin{bmatrix} x' \\ t' \end{bmatrix}=\begin{bmatrix} \cos\alpha & \sin\alpha \\ -\sin\alpha & \cos\alpha \end{bmatrix}\begin{bmatrix} x \\ t \end{bmatrix}.$$

This relation describes a rotation of the vector $(x,t)$ on the $xt$-plane by an angle $\alpha$. The PDE $F[u]=0$ is thus invariant under such a rotation on the plane of the independent variables.

## 10. Bäcklund-transformation-related recursion operators

As defined in Sec. 6, a recursion operator is a linear operator that produces a new symmetry characteristic of a PDE when it acts on an "old" characteristic. As seen in the example of Sec. 7, a recursion operator can also be viewed as an auto-Bäcklund transformation (auto-BT) for the symmetry condition of the PDE, which condition is a linear PDE for the symmetry characteristic, valid for solutions of the original PDE.

It is known that when two (generally nonlinear) PDEs are connected by a BT, symmetries of one PDE can be associated with symmetries of the other. In [9] we went one step further by asking the question: can we find a transformation that relates *recursion operators* of these PDEs? This would be a relation connecting two auto-BTs, each of which relates symmetries of a respective PDE. In other words, we were after "a transformation of transformations" rather than a transformation of functions.

Our "laboratory" model was the *self-dual Yang-Mills equation* (SDYM) and its potential version, the *potential SDYM equation* (PSDYM). The two PDEs are related by a BT and it was shown in [15] that the symmetries of PSDYM can be used to construct *potential symmetries* of SDYM. As proven in [9], SDYM and PSDYM have isomorphic infinite-dimensional symmetry Lie subalgebras generated by BT-related recursion operators.

The SDYM equation is written

$$F[J] \equiv (J^{-1}J_y)_{\bar{y}} + (J^{-1}J_z)_{\bar{z}} = 0 \tag{70}$$

where the $y, z, \bar{y}, \bar{z}$ are complex variables which will be collectively denoted $x^\mu$ ($\mu$=1,...,4). As usual, subscripts indicate total derivatives $D_\mu$ with respect to the $x^\mu$. We assume that the dependent variable $J$ is $SL(N,C)$-valued, i.e., is represented by a complex ($N\times N$) matrix with $\det J=1$.

We consider a system of PDEs involving $J$ and another matrix function $X$ that plays the role of the "potential" of (70):

$$J^{-1}J_y = X_{\bar{z}} \quad , \quad J^{-1}J_z = -X_{\bar{y}} \tag{71}$$

The system (71) is a BT relating $J$ with $X$. The integrability condition $(X_{\bar{y}})_{\bar{z}} = (X_{\bar{z}})_{\bar{y}}$ yields the SDYM equation (70). Another integrability condition is found by applying the matrix identity (45):

$$(J^{-1}J_z)_y - (J^{-1}J_y)_z + [J^{-1}J_y \, , \, J^{-1}J_z] = 0 \, ,$$

which yields the PSDYM equation for $X$:

$$G[X] \equiv X_{y\bar{y}} + X_{z\bar{z}} - [X_{\bar{y}} \, , \, X_{\bar{z}}] = 0 \tag{72}$$

As follows from (71), the condition $\det J=1$ can be satisfied by requiring that $trX$=0, which is compatible with (72). Thus, $SL(N,C)$ SDYM solutions correspond to $sl(N,C)$ PSDYM solutions through the BT (71).

For this problem we need an enhanced jet space with variables $x^\mu$, $J$, $X$ and partial derivatives of various orders of $J$ and $X$ with respect to the $x^\mu$. Functions $M[J;X]$ in this space are assumed to be of a form generalizing (1):

$$M[J;X] = \sum a(x^\mu)\,\Pi[J;X]\,b(x^\mu)\,\Pi'[J;X]\,c(x^\mu)\cdots \tag{73}$$

where the matrices $\Pi[J;X]$, $\Pi'[J;X]$, etc., are products of powers and/or derivatives of $J$ and $X$ with respect to the $x^\mu$. Characteristic derivatives acting on functions (73) are of the form

$$\Delta = \Delta_Q + \Delta_\Phi \tag{74}$$

where

$$\Delta_Q J = Q[J;X] \ , \quad \Delta_\Phi X = \Phi[J;X] \quad \text{and} \quad \Delta_Q X = \Delta_\Phi J = 0 \ .$$

Therefore,

$$\Delta J = Q[J;X] \ , \quad \Delta X = \Phi[J;X] \tag{75}$$

for given characteristic functions $Q[J;X]$ and $\Phi[J;X]$. Clearly, the operator $\Delta$ has all the properties of characteristic derivatives stated in Sec. 2. For example,

$$\begin{gathered}
\Delta f(x^\mu) = 0 \ , \\
\Delta D_\mu = D_\mu \Delta \ , \\
\Delta(MN) = (\Delta M)N + M(\Delta N) \ , \\
\Delta[M,N] = [\Delta M, N] + [M, \Delta N] \ , \\
\Delta(M^{-1}) = -M^{-1}(\Delta M)M^{-1} \ , \text{ etc.}
\end{gathered}$$

We introduce the "covariant derivative" operators (a term borrowed from gauge theories):

$$\begin{aligned}
\hat{A}_y &= D_y + [J^{-1}J_y \, , \ \ ] = D_y + [X_{\bar{z}} \, , \ \ ] \\
\hat{A}_z &= D_z + [J^{-1}J_z \, , \ \ ] = D_z - [X_{\bar{y}} \, , \ \ ]
\end{aligned} \tag{76}$$

where the BT (71) has been taken into account. As can be shown, the "zero curvature" condition $[\hat{A}_y , \hat{A}_z] = 0$ is satisfied, consistently with the fact that the "connections" $J^{-1}J_y$ and $J^{-1}J_z$ are pure gauges. The linear operators (76) are derivations on the Lie algebra of $sl(N,C)$-valued functions, satisfying the Leibniz rule in the form

$$\begin{gathered}
\hat{A}_y [M,N] = [\hat{A}_y M, N] + [M, \hat{A}_y N] \ , \\
\hat{A}_z [M,N] = [\hat{A}_z M, N] + [M, \hat{A}_z N] \ .
\end{gathered}$$

Let now $\delta J = \alpha Q[J;X]$ and $\delta X = \alpha \Phi[J;X]$ represent an infinitesimal symmetry of the system (71), where $\alpha$ is an infinitesimal parameter. This means that the pair of functions $(J+\delta J, X+\delta X)$ is a solution of the system if $(J, X)$ is a solution. This, in turn, suggests that the integrability conditions $F[J+\delta J]=0$ and $G[X+\delta X]=0$ are satisfied when the integrability conditions $F[J]=0$ and $G[X]=0$ are satisfied; that is, $J+\delta J$ and $X+\delta X$ are solutions of (70) and (72), respectively, when $J$ and $X$ are solutions of these PDEs. The functions $Q[J;X]$ and $\Phi[J;X]$ are thus, respectively, SDYM and PSDYM symmetry characteristics. The characteristic derivative $\Delta$ with respect to these functions is given by Eqs. (74) and (75) and represents a symmetry of the BT (71), thus also of the PDEs (70) and (72).

The symmetry condition for the SDYM equation (70) is $\Delta F[J]=0$ (mod $F[J]$, a condition that will not be stated explicitly in the sequel), which yields

$$D_{\bar{y}}\Delta(J^{-1}J_y)+D_{\bar{z}}\Delta(J^{-1}J_z)=0 \ .$$

By using (76) and the fact that $\Delta J=Q$ it can be shown that

$$\Delta(J^{-1}J_y)=\hat{A}_y(J^{-1}Q) \ , \quad \Delta(J^{-1}J_z)=\hat{A}_z(J^{-1}Q) \qquad (77)$$

The SDYM symmetry condition then becomes

$$(D_{\bar{y}}\hat{A}_y+D_{\bar{z}}\hat{A}_z)(J^{-1}Q)=0 \qquad (78)$$

The symmetry condition for PSDYM (72) is $\Delta G[X]=0$, which yields, by using (76) and the fact that $\Delta X=\Phi$:

$$\hat{A}_y\Phi_{\bar{y}}+\hat{A}_z\Phi_{\bar{z}}\equiv(\hat{A}_y D_{\bar{y}}+\hat{A}_z D_{\bar{z}})\Phi=0 \qquad (79)$$

The following operator identity can be proven:

$$\hat{A}_y D_{\bar{y}}+\hat{A}_z D_{\bar{z}}=D_{\bar{y}}\hat{A}_y+D_{\bar{z}}\hat{A}_z$$

by which (79) is written in the alternative form

$$(D_{\bar{y}}\hat{A}_y+D_{\bar{z}}\hat{A}_z)\Phi=0 \qquad (80)$$

By comparing (78) and (80) we observe that the PSDYM symmetry characteristic $\Phi$ and the function $J^{-1}Q$, where $Q$ is an SDYM characteristic, satisfy the same symmetry condition. We thus conclude the following:

- If $Q$ is an SDYM characteristic, then $\Phi=J^{-1}Q$ is a PSDYM characteristic.

Conversely,

- if $\Phi$ is a PSDYM characteristic, then $Q=J\Phi$ is an SDYM characteristic.

Finally, by using (75) and (77) we find a pair of PDEs that express the symmetry condition for the BT (71):

$$\hat{A}_y(J^{-1}Q)=\Phi_{\bar{z}} \ , \quad \hat{A}_z(J^{-1}Q)=-\Phi_{\bar{y}} \qquad (81)$$

The system (81) is a BT relating a symmetry characteristic $\Phi$ of PSDYM with a symmetry characteristic $Q$ of SDYM. Indeed, the integrability condition $(\Phi_{\bar{z}})_{\bar{y}}=(\Phi_{\bar{y}})_{\bar{z}}$ yields the symmetry condition (78) of SDYM, while the integrability condition $[\hat{A}_y,\hat{A}_z](J^{-1}Q)=0$, valid in view of the commutativity of the covariant derivatives (76), yields the PSDYM symmetry condition (79).

Now, as mentioned previously, the product $J^{-1}Q$ is a PSDYM characteristic, which we call $\Phi_0$. We thus rewrite the first equation in the system (81) as

$$\Phi_{\bar{z}} = \hat{A}_y \Phi_0 \quad ,$$

which is formally integrated to give

$$\Phi = D_{\bar{z}}^{-1} \hat{A}_y \Phi_0 \equiv \hat{R} \Phi_0 \tag{82}$$

where we have introduced the linear operator

$$\hat{R} = D_{\bar{z}}^{-1} \hat{A}_y \tag{83}$$

As suggested by (82), the operator (83) is a *recursion operator* for PSDYM, producing a new symmetry characteristic $\Phi$ when acting on an "old" characteristic $\Phi_0$.

Let $\Phi$ and $\Phi_0$ be PSDYM symmetry characteristics related to each other by (82). Then $Q=J\Phi$ and $Q_0=J\Phi_0$ are SDYM characteristics, in terms of which (82) is rewritten as

$$Q = J \hat{R} J^{-1} Q_0 \equiv \hat{T} Q_0 \tag{84}$$

where the linear operator

$$\hat{T} = J \hat{R} J^{-1} \tag{85}$$

is a recursion operator for SDYM. This operator was constructed previously by using a Lax pair for SDYM [11].

Equation (85) can be regarded as a transformation relating recursion operators of two PDEs – namely, (70) and (72) – connected to each other by the BT (71). As shown in [9], the recursion operators (83) and (85) generate isomorphic, infinite-dimensional symmetry Lie algebras of PSDYM and SDYM, respectively. In practical terms, infinite-dimensional symmetry subalgebras of SDYM, such as Kac-Moody and Virasoro algebras, can be obtained directly from the easier-to-find corresponding symmetries of PSDYM [20].

## 11. Concluding remarks

The algebraic approach to the symmetry problem of PDEs, presented in this article, is, in a sense, an extension to matrix-valued problems of the ideas contained in [1], in much the same way as [4] and [5] constitute a generalization of the Harrison-Estabrook geometrical approach [2] to matrix-valued (as well as vector-valued and Lie-algebra-valued) PDEs.

The symmetry transformations we have considered involve only the change of the dependent variable of the PDE while leaving the independent variables unchanged. Indeed, as Olver [1] has shown, *every* symmetry of a PDE may be expressed as a transformation of the dependent variable alone, i.e., as a transformation in the fiber space. The symmetry-generating characteristic derivative $\Delta_Q$ corresponds to Olver's *evolutionary vector field* with characteristic $Q$.

In *local* (point) symmetries of the PDE $F[u]=0$ the symmetry characteristics $Q[u]$ contain at most first-order derivatives of $u$ with respect to the independent variables (a number of such symmetries were considered in Sec. 9 in connection with the Laplace equation). The case of *generalized* (non-local) symmetries is more complex; a formal solution to the problem of obtaining one-parameter families of generalized symmetry transformations of PDEs is given in Sec. 5.1 of [1].

Admittedly, the abstract algebraic formalism we have presented does not exhibit significant advantages over the standard geometrical methods (in particular, those described in [4]) with regard to *finding* symmetries of PDEs. However, by employing the concept of the characteristic derivative one is able to bypass the difficulty of having to represent symmetry-generating operators as vector fields in the form of differential operators when matrix-valued variables are involved, which situation can only be handled by making certain *ad hoc* assumptions regarding the action of such "unnatural" operators. The algebraic view offers a more rigorous framework for identifying symmetry operators and finding infinitesimal symmetry transformations of matrix-valued PDEs, as well as for studying the Lie-algebraic structure of the set of symmetry generators [9,20].